\documentclass[3p,times]{elsarticle}
\usepackage{amsmath}
\usepackage{mathrsfs}
\usepackage{graphicx}
\usepackage{amssymb}
\usepackage[figuresright]{rotating}
\usepackage{booktabs, multirow, caption}
\usepackage{multicol}
\usepackage{float}
\usepackage{subfigure}
\usepackage{color}
\usepackage{epstopdf}
\usepackage{makecell}
\usepackage{amsthm}                     
\theoremstyle{elsarticle-bold}

\newtheorem{lemma}{Lemma}[section]

\newtheorem{theorem}{Theorem}[section]

\newtheorem{remark}{Remark}[section]

\newtheorem{example}{Example}[section]

\begin{document}

\begin{frontmatter}

\title{The double splitting iteration method\\ for solving the large indefinite least squares problem}

\author{Jun Li$^{a}$\footnote{Corresponding author. Email: junli026430@163.com.}, Lingsheng Meng$^{b}$}
\address{$^a$~School of Science, Lanzhou University of Technology, Lanzhou, Gansu, 730050, P. R. China\\
$^b$~College of Mathematics and Statistics, Northwest Normal University, Lanzhou, Gansu, 730070, P. R. China}

\begin{abstract}
Addressing large-scale indefinite least squares (ILS) problem poses notable computational bottlenecks in the field of numerical linear algebra. State-of-the-art iterative schemes for such problems are predominantly constructed upon the single splitting of the coefficient matrix derived from the corresponding normal equation. In this work, we put forward an innovative iterative framework grounded in the double splitting of normal equations tailored for ILS problem. Specifically, we elaborate on a distinct implementations of the double splitting strategy, which offer constructive insights and methodological references for subsequent research on double splitting-based iterative methods. Two numerical experiments further corroborate that the proposed double splitting iterative paradigm outperforms conventional single splitting approaches in both computational efficiency and convergence robustness.

\end{abstract}


\end{frontmatter}

\section{Introduction}
\pagestyle{plain} \setcounter{equation}{ 0}
\renewcommand{\theequation}{1.\arabic{equation}}

The indefinite least squares (ILS) problem arises in a variety of practical applications, including total least squares \cite{Golub1980,Huffel1991,Samar2021,Samar2022} and $H^\infty$ smoothing in optimization \cite{Hassibi1993,Sayed1996}. This class of problem was first introduced in \cite{Chan1998} and has since attracted considerable research attention \cite{Bojan2021,Diao2019,Li2018,Liu2011,Liu2013,Liu2014,Xu2004}. In this work, we focus on solving the ILS problem
\begin{equation} \label{1.1}
	\mathrm{ILS:}~\min_{x\in \mathbb{R}^{n}}(b-Ax)^TJ(b-Ax),
\end{equation}
where $A\in \mathbb{R}^{m\times n}$ is a large sparse matrix with $m\geq n$, $b\in \mathbb{R}^{m}$, and $J=\operatorname{diag}(I_p,-I_q)$ denotes the signature matrix with $p+q=m$. When either $p=0$ or $q=0$, problem (\ref{1.1}) reduces to the standard least squares (LS) problem, whose quadratic form is definite. For the case $pq>0$, the problem involves minimizing an indefinite quadratic form associated with $J$, leading to the normal equation
\begin{equation}\label{1.2}
	A^TJAx = A^TJb.
\end{equation}
The Hessian matrix of the ILS problem is $2A^TJA$, and the uniqueness of the solution is determined by whether $A^TJA$ is symmetric positive definite (SPD).

For small and dense ILS problems, direct methods based on QR decomposition, such as Cholesky factorization or Householder transformations, have been employed \cite{Bojan2003,Bojan2021,Chan1998,Xu2004}. In 2011, Liu et al. \cite{Liu2011} developed a preconditioner based on a partition of $A$ and applied preconditioned conjugate gradient methods. Subsequently, they introduced incomplete hyperbolic Gram–Schmidt-based preconditioners \cite{Liu2013} by approximating $A^TJA$ via incomplete hyperbolic classical or modified Gram–Schmidt procedures. In 2014, a block successive over-relaxation (SOR) iteration method with a relaxation parameter was proposed \cite{Liu2014}. To further accelerate the convergence of the SOR method, Song et al. \cite{Song2020} proposed the USSOR iteration method for solving the block $3\times 3$ linear systems arising from the ILS problem. 

For large and sparse ILS problem, iteration methods are often more attractive. Among iterative solvers, Krylov subspace methods (such as the generalized minimal residual method, GMRES) have gained widespread attention due to their excellent convergence properties, especially when dealing with large-scale sparse linear systems. However, their numerical performance can be significantly enhanced—rather than merely better demonstrated—when paired with a suitably designed preconditioner. A good preconditioner transforms the original linear system into one with more favorable spectral properties, thereby accelerating the convergence rate of Krylov subspace methods and reducing the overall computational cost. In recent years, many researchers have proposed new preconditioners within the GMRES framework to solve three-by-three block linear systems, and thus to obtain the solution of the ILS problem (\ref{1.1}); see \cite{Hashemi2026,Khojasteh2026,Li2025,Li2026,Xin2025} for details. In addition to Krylov subspace methods, the stationary iteration method has also been studied. In 2023, Zhang and Li \cite{Zhang2023} presented splitting-based (SP) randomized iteration methods for normal equation (\ref{1.2}). More recently, in 2025, Meng et al. \cite{Meng2025a} introduced a variable Uzawa iteration method for three-by-three block linear system arising from ILS problem, and in the same year they also proposed a generalized SP (GSP) iteration method and an alternating direction implicit (ADI) iteration method for the equation (\ref{1.2}) \cite{Meng2025b}. The list of relevant studies is extensive and cannot be fully enumerated here.

For existing iterative solvers for large ILS problem, the coefficient matrix of the normal equation (\ref{1.2}) is first split as
$$
A^TJA = M - N,
$$
where $M$ is nonsingular. This leads to the standard splitting iteration scheme
\begin{align}\label{1.3}
	Mx^{k+1} = Nx^k + A^TJb,\qquad k=0,1,2,\dots
\end{align}
In particular, the matrix $A$ is often partitioned as $A = [A_1;\,A_2]$ with $A_1 \in \mathbb{R}^{p\times n}$ of full column rank and $A_2 \in \mathbb{R}^{q\times n}$. Consequently, $A_1^TA_1$ is SPD and $A_2^TA_2$ is symmetric positive semidefinite. Based on this partition, several effective splittings have been proposed:

\begin{itemize}
	\item When $M = A_1^TA_1$ and $N = A_2^TA_2$, scheme (\ref{1.3}) is called the splitting-based (SP) iteration method, which is unconditionally convergent \cite{Zhang2023}.
	\item When $M = \alpha I + A_1^TA_1$ and $N = \alpha I + A_2^TA_2$ with $\alpha>0$, it is referred to as the generalized SP (GSP) iteration method, and convergence holds for all $\alpha>0$ \cite{Meng2025b}.
	\item A more elaborate alternating scheme arises by taking
	$$
	M_1 = \alpha I + A_1^TA_1,\quad N_1 = \alpha I + A_2^TA_2,\quad
	M_2 = -(\beta I + A_2^TA_2),\quad N_2 = -(\beta I + A_1^TA_1)
	$$
	with $\alpha,\beta>0$. This yields the alternating direction implicit (ADI) iteration method:
	$$
	\begin{cases}
		M_1 x^{k+\frac{1}{2}} = N_1 x^k + A^TJb,\\[2mm]
		M_2 x^{k+1} = N_2 x^{k+\frac{1}{2}} + A^TJb.
	\end{cases}
	$$
The ADI method converges under the condition $\beta > \alpha > 0$ \cite{Meng2025b}.
\end{itemize}
The above three splittings are all instances of \emph{single splittings}. Accordingly, the SP and GSP methods are one-step linear stationary iteration methods, while the ADI method is a two-half-step method. Numerical experiments in the literature confirm their effectiveness for solving large ILS problem.

Inspired by the work of Wo$\acute{\mbox{z}}$nicki \cite{Woznicki1993}, we consider a different decomposition of the coefficient matrix:
$$
A^TJA = P - R - S,
$$
where $P$ is nonsingular. This is called a \emph{double splitting} of $A^TJA$. Based on this decomposition, we derive a two-step iteration scheme of the form
\begin{align*}
	x^{k+1} = P^{-1}R x^k + P^{-1}S x^{k-1} + P^{-1}A^TJb,\qquad k=1,2,\dots
\end{align*}
Such two-step methods have the potential to achieve faster convergence than standard one-step methods, as they incorporate information from two previous iterates.

In this paper, we propose a concrete double splitting of $A^TJA$. It's convergence properties are analyzed in Section~2. Section~3 presents numerical results that demonstrate the superiority of the proposed double splitting iteration method over the single splitting methods (SP, GSP, and ADI) in terms of convergence speed. Finally, concluding remarks are given in Section~4.

In the end, we introduce some used notations in the paper. Matlab notation $[x;y;z]$ represents the column vector $[x^T,y^T,z^T]^T$. For a given nonsingular matrix $Z$, $\rho(Z)$ denotes the spectral radius of matrix $Z$, $Z\succ 0$ means that $Z$ is symmetric positive definite (SPD) matrix. The conjugate transpose of a vector $\ell$ is expressed as $\ell^*$. A matrix or vector whose elements are all zeros is defined as $\mathbf{0}$.

\section{The double splitting iteration method}
\pagestyle{plain} \setcounter{equation}{ 0}
\renewcommand{\theequation}{2.\arabic{equation}}

In this section, we propose a double splitting strategy for the normal equation (\ref{1.2}), which serves as the foundation for a two-step iteration method to obtain the solution of the ILS problem (\ref{1.1}).

Following the concept of double splitting, we split the matrix $A^TJA$ as
\begin{align}\label{2.1}
	A^TJA := P - R - S,
\end{align}
where $P$ is required to be nonsingular. Inspired by the splitting pattern of the GSP method, we choose the following specific matrices:
\begin{align}\label{2.2}
	P = \alpha I + A_1^TA_1,\qquad R = A_2^TA_2,\qquad S = \alpha I,
\end{align}
with $\alpha > 0$ a given parameter. In this setting, $P$ is SPD, while $R$ and $S$ are symmetric. Based on the double splitting (\ref{2.1}), the solution of the ILS problem can be obtained via the following two-step stationary iteration:
\begin{align}\label{2.3}
	x^{k+1} = P^{-1}R x^{k} + P^{-1}S x^{k-1} + P^{-1}A^TJb,\qquad k=1,2,\dots
\end{align}

To analyze the convergence of scheme (\ref{2.3}), we rewrite it in an equivalent block form. Introducing the augmented vectors $\mathbf{x}^k = [ x^{k+1}; x^{k}]$, we obtain
\begin{align}\label{2.4}
	\begin{bmatrix} x^{k+1} \\ x^{k} \end{bmatrix}
	=
	\begin{bmatrix} P^{-1}R & P^{-1}S \\ I & \mathbf{0} \end{bmatrix}
	\begin{bmatrix} x^{k} \\ x^{k-1} \end{bmatrix}
	+
	\begin{bmatrix} P^{-1}A^TJb \\ \mathbf{0} \end{bmatrix},\qquad k=1,2,\dots
\end{align}
We refer to this method as the \emph{double splitting (DS) iteration method}. Denote the iteration matrix by
\[
W = \begin{bmatrix} P^{-1}R & P^{-1}S \\ I & \mathbf{0} \end{bmatrix}.
\]
It is well known  that for any initial vectors $x^0, x^1$, the sequence $x^{k+1}$ generated by (\ref{2.4}) converges to the unique solution of (\ref{1.1}) if and only if the spectral radius of $W$ satisfies $\rho(W) < 1$ (see, e.g., \cite{Benzi2005,Saad2003}). In this case, the double splitting (\ref{2.1}) is said to be convergent.

To establish conditions under which $\rho(W) < 1$ holds for the DS iteration method, we study the eigenvalue distribution of $W$. Before proceeding to the main convergence analysis, we first recall a useful lemma.

\begin{lemma}\cite{Young1971}\label{lem2.1}
Both roots of the real quadratic equation $x^2-px+q=0$ are less than $1$ in modulus if and only if $\mid q \mid < 1$ and $\mid p\mid < 1+q$.
\end{lemma}

Let $(\lambda, \ell)$ with $\ell = (\ell_1; \ell_2)$ be an eigenpair of the iteration matrix $W$, i.e., $W\ell = \lambda \ell$. Using the block representation
\[
W = \begin{bmatrix} P^{-1}R & P^{-1}S \\ I & \mathbf{0} \end{bmatrix}
= \begin{bmatrix} P & \mathbf{0} \\ \mathbf{0} & I \end{bmatrix}^{-1}
\begin{bmatrix} R & S \\ I & \mathbf{0} \end{bmatrix},
\]
the eigenvalue equation becomes
\[
\begin{bmatrix} R & S \\ I & \mathbf{0} \end{bmatrix}
\begin{bmatrix} \ell_1 \\ \ell_2 \end{bmatrix}
= \lambda \begin{bmatrix} P & \mathbf{0} \\ \mathbf{0} & I \end{bmatrix}
\begin{bmatrix} \ell_1 \\ \ell_2 \end{bmatrix},
\]
which is equivalent to the system
\begin{equation}\label{2.5}
	\left\{
	\begin{array}{ll}
		R\ell_1 + S\ell_2 = \lambda P\ell_1,\\[2mm]
		\ell_1 = \lambda \ell_2.
	\end{array}
	\right.
\end{equation}

If $\lambda = 1$, then $\ell_1 = \ell_2$. Substituting into the first equation of (\ref{2.5}) yields
$(P - R - S)\ell_1 = A^TJA\ell_1 = \mathbf{0}$. Since $A^TJA$ is SPD (as the ILS solution is unique), we obtain $\ell_1 = \mathbf{0}$, and consequently $\ell_2 = \mathbf{0}$, contradicting that $\ell$ is an eigenvector. Hence $\lambda \neq 1$.

When $\lambda \neq 1$, we must have $\ell_2 \neq \mathbf{0}$; otherwise $\ell$ would be the zero vector. Substituting $\ell_1 = \lambda \ell_2$ into the first equation in (\ref{2.5}) gives
\begin{align}\label{2.6}
	\left(\lambda^2 P - \lambda R - S\right)\ell_2 = \mathbf{0}.
\end{align}
Define the following Rayleigh quotients with respect to $\ell_2 \neq \mathbf{0}$:
\[
a := \frac{\ell_2^* P \ell_2}{\ell_2^* \ell_2}
= \frac{\ell_2^*(\alpha I + A_1^TA_1)\ell_2}{\ell_2^* \ell_2} > 0,
\quad
a_2 := \frac{\ell_2^* R \ell_2}{\ell_2^* \ell_2}
= \frac{\ell_2^* A_2^TA_2 \ell_2}{\ell_2^* \ell_2} \ge 0,
\quad
a_3 := \frac{\ell_2^* S \ell_2}{\ell_2^* \ell_2} = \alpha > 0.
\]
Multiplying (\ref{2.6}) on the left by $\frac{\ell_2^*}{\ell_2^*\ell_2}$ and using the Rayleigh quotients yields the scalar quadratic equation
\begin{align}\label{2.7}
	a \lambda^2 - a_2 \lambda - a_3 = 0.
\end{align}

Recall that $a = \alpha + a_1$ with $a_1 = \frac{\ell_2^* A_1^TA_1 \ell_2}{\ell_2^* \ell_2} > 0$, and $a_3 = \alpha$. Then (\ref{2.7}) becomes
\[
(\alpha + a_1)\lambda^2 - a_2 \lambda - \alpha = 0,
\]
or equivalently
\[
\lambda^2 - \frac{a_2}{\alpha + a_1}\,\lambda - \frac{\alpha}{\alpha + a_1} = 0.
\]

To ensure $\rho(W) < 1$, we need all eigenvalues $\lambda$ to satisfy $|\lambda| < 1$. Lemma \ref{lem2.1} shows that $|\lambda| < 1$ if and only if $|c| < 1$ and $|b| < 1 + c$, here $b = \frac{a_2}{\alpha + a_1}$ and $c =- \frac{\alpha}{\alpha + a_1}$. The condition $|c| < 1$ holds trivially because $\frac{\alpha}{\alpha + a_1} < 1$ for any $\alpha > 0$ and $a_1 > 0$. Moreover,
\[
|b| = \frac{a_2}{\alpha + a_1} < \frac{a_1}{\alpha + a_1} = 1 + \frac{-\alpha}{\alpha + a_1}.
\]
Therefore, both conditions are satisfied for all $\alpha > 0$, implying that every eigenvalue $\lambda$ of $W$ satisfies $|\lambda| < 1$. Consequently, the DS iteration scheme (\ref{2.3}) or (\ref{2.4}) converges unconditionally with respect to the parameter $\alpha$.

The above result is summarized in the following theorem.

\begin{theorem}\label{thm2.1}
	Consider the DS iteration scheme (\ref{2.3}) (or equivalently (\ref{2.4})) for solving the normal equation (\ref{1.2}), which yields the solution of the ILS problem (\ref{1.1}). For any positive parameter $\alpha > 0$, the DS iteration method is convergent.
\end{theorem}

\begin{remark}
	Since Wo$\acute{\mbox{z}}$nicki \cite{Woznicki1993} introduced the concept of double splitting of matrices, many researchers have investigated double splittings for various matrix classes, such as Hermitian positive definite matrices \cite{Shen2007} and non‑Hermitian positive semidefinite matrices \cite{Zhang2010}. These works also discussed the advantages of double splittings over single splittings. In particular, \cite[Theorem 2.1]{Shen2007} shows that if both $P + S$ and $A^TJA + 2R$ are SPD, then the two-step iteration derived from the double splitting is convergent. For our choice (\ref{2.2}), it is easy to verify that $P+S = \alpha I + A_1^TA_1 + \alpha I = 2\alpha I + A_1^TA_1$ is SPD, and so is $A^TJA+2R = A_1^TA_1 - A_2^TA_2 + 2A_2^TA_2 = A_1^TA_1 + A_2^TA_2$. Hence, the convergence of the DS iteration method can also be deduced from \cite{Shen2007}. Nevertheless, our direct spectral analysis confirms that the method converges for all $\alpha>0$.
\end{remark}

Based on the double splitting defined in (\ref{2.2}), the algorithmic implementation of the DS iteration scheme is described as follows:
\begin{table}[!ht]
\renewcommand\arraystretch{1.2}
\begin{tabular*}{\textwidth}{@{\extracolsep{\fill}}lllllllllllllllllllllllllllllllllllll} 
\Xhline{1.5pt}
Algorithm 2.1. The DS iteration method for ILS problem (\ref{1.1}).\\	
\hline
1. Input $A$, $J$, $b$ and initial iteration vectors $x^0$ and $x^1$:\\
2. Calculate $Q_1=(\alpha I+A_1^TA_1)^{-1}$ by Cholesky decomposition;\\
3. Set $\bar{R}_1=Q_1A_2^TA_2,$ $\bar{S}_1=\alpha Q_1$, $B_1=Q_1A^TJb$;\\
4. $\mathbf{for}$ $k=1,2,\cdots$, until converges,  $\mathbf{do}$\\
5. Update $x^{k+1}=\bar{R}_1x^{k}+\bar{S}_1x^{k-1}+B_1$;\\
6.  $\mathbf{end~for}$.\\
\Xhline{1.5pt}
\end{tabular*}
\end{table}

\section{Numerical experiments}

In this section, we present two numerical examples to illustrate the effectiveness of the proposed double splitting (DS) iteration method compared with existing single splitting methods, namely the SP, GSP and ADI iteration methods. The formulations of these single splitting methods are described in Section~1, and their detailed algorithmic implementations can be found in \cite{Zhang2023,Meng2025b}. In particular, we select $\alpha=10^{-6}$ for GSP method, $\alpha=10^{-6}$ ans $\beta=10^{15}$ are applied to ADI method. We solve the normal equation $A^TJAx = A^TJb$ to obtain the solution of the ILS problem (\ref{1.1}) for a fair comparison. 

For the SP, GSP and ADI methods, the initial iteration vector is set to the zero vector. For the DS method, we take two zero vectors as the initial guesses, i.e., $x^0 = x^1 = \mathbf{0}$. All iterations are terminated when the relative residual satisfies
\[
\text{RES} = \frac{\|A^TJb - A^TJA x^{k}\|_2^2}{\|A^TJb\|_2^2} < 10^{-8},
\]
where $x^{k}$ denotes the approximate solution at the $k$th iteration, or when the maximum number of iterations $k_{\max}=10000$ is reached. 

The numerical performance of the four methods is evaluated in terms of the number of iteration steps (denoted as `IT'), the computational time in seconds (denoted as `CPU'). All experiments are carried out using MATLAB (R2018a) on a personal computer equipped with an AMD i5-13400F CPU at 2.5 GHz and 32.00 GB of RAM.

\begin{figure}[htbp]
	\centering
	\includegraphics[height=3.5cm,width=4cm]{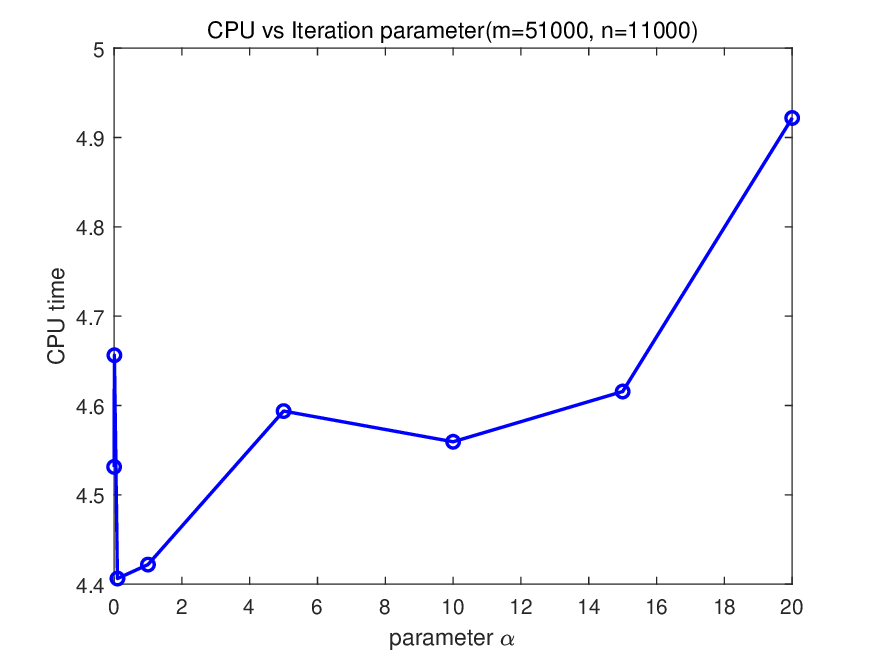}
	\includegraphics[height=3.5cm,width=4cm]{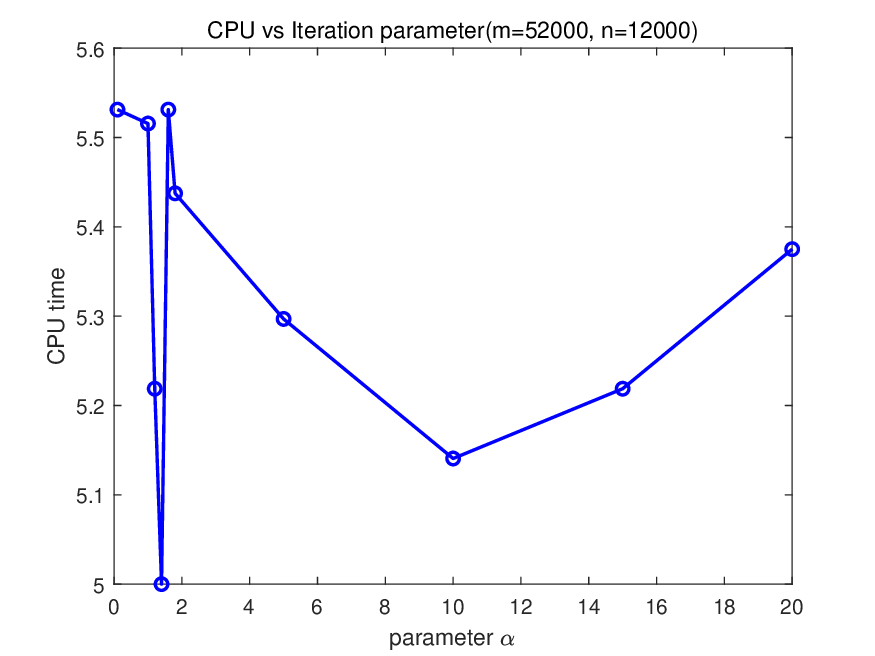}
	\includegraphics[height=3.5cm,width=4cm]{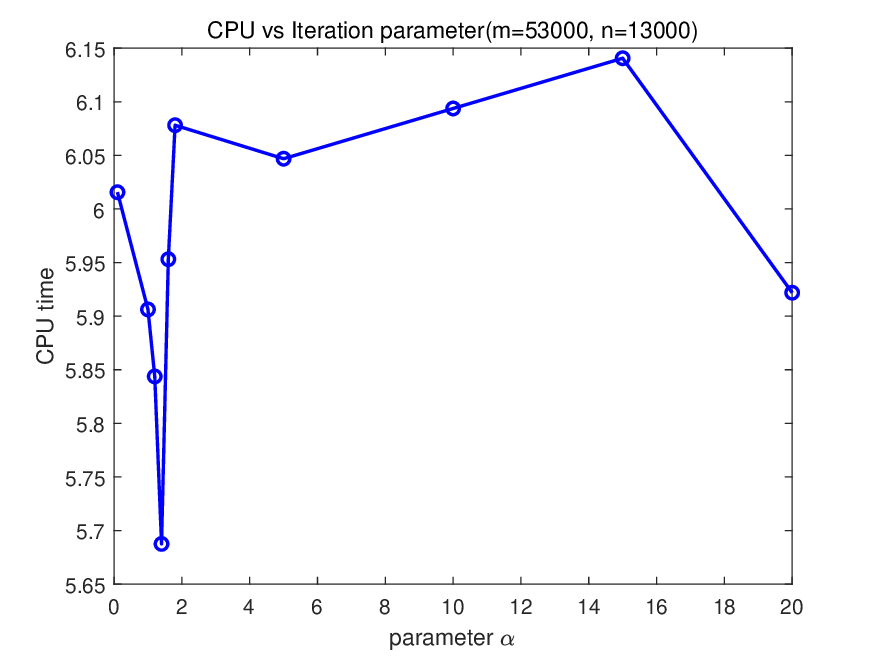}
	\includegraphics[height=3.5cm,width=4cm]{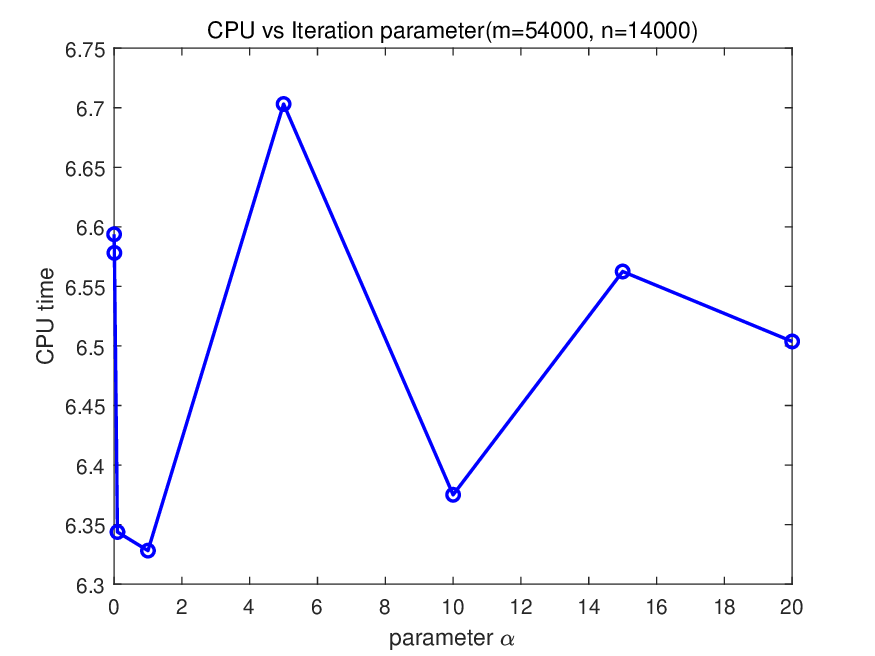}
	\caption{\footnotesize \emph{$\alpha$ vs CPU times of DS iteration method in Example 4.1.}}\label{fig1}
\end{figure}

\begin{example}\cite{Meng2025b}
we construct the ILS problem (\ref{1.1}) using random large scale dense matrix $A_1 = rand(p, n)$, and $A_2 =
7*eye(q, n)$, $b_1 = rand(p, 1)$, $b_2 = rand(q, 1)$, where $p = 40000$, $n=q = [11000 : 1000 : 14000]$ and $m=p+q$.
\end{example}

\begin{table}[htbp]  
	\centering
	\captionsetup{width=0.9\textwidth, font=small, labelfont=bf, labelsep=period}
	\renewcommand{\arraystretch}{1}  
	\setlength{\tabcolsep}{5pt}       
	\caption{The numerical results of various iteration methods with dense matrices for Example 4.1.}\label{tab3}
	
	\begin{tabular}{l *{8}{r}}
		\toprule
		\multirow{2}{*}{Method} 
		& \multicolumn{2}{c}{$m=51000,\; n=11000$} 
		& \multicolumn{2}{c}{$m=52000,\; n=12000$} 
		& \multicolumn{2}{c}{$m=53000,\; n=13000$}
		& \multicolumn{2}{c}{$m=54000,\; n=14000$} \\
		\cmidrule(lr){2-3} \cmidrule(lr){4-5} \cmidrule(lr){6-7} \cmidrule(lr){8-9}
		& IT & CPU (s) & IT & CPU (s) & IT & CPU (s) & IT & CPU (s) \\
		\midrule
		SP   & 1 & 77.2813 & 1 & 97.0156 & 1 & 121.5625 & 1 & 140.3125 \\
		GSP  & 1 & 75.1406 & 1 & 96.8594 & 1 & 118.0938 & 1 & 139.5781 \\
		ADI  & 1 & 13.6094 & 1 & 17.6406 & 1 & 24.8438  & 1 & 29.9375\\
		DS  & 2 &  4.8125 & 2 &  5.1563 & 2 &  5.5469   & 2 &  6.4531 \\
		\bottomrule
	\end{tabular}
\end{table}

Figure~1 illustrates the relationship between the choice of the parameter $\alpha$ and the CPU time of the DS iteration method under different matrix dimensions. Based on the observed performance in this figure, we select $\alpha = 1$ as the optimal parameter for the DS iteration method in Example~3.1 (i.e., for the numerical tests reported therein). Table~1 presents a comparative evaluation of different iteration methods for solving the ILS problem (\ref{1.1}). The results clearly demonstrate that, with the chosen optimal parameter, the proposed double-splitting iteration method converges in only two iterations. Moreover, its solution time is significantly lower than that of the other iteration methods considered, highlighting the computational efficiency of the DS approach.

\begin{example}\cite{Li2026,Liu2011,Xin2025}
Let $p,~q,~ n$ take different values, $\epsilon$ is given parameter and $
	\tilde{B}=Y\cdot\mbox{diag}(D,\mathbf{0})\cdot Z^{T}\in  \mathbb{R}^{p\times n}$, 
where $Y\in  \mathbb{R}^{p\times p}$, $Z\in  \mathbb{R}^{n\times n}$ are given orthogonal matrices and $ D={\rm diag}(1,\frac{1}{2},..., \frac{1}{n})\in  \mathbb{R}^{n\times n}$. Next let $B=\tilde{B}+\epsilon E$, $d=\tilde{B}1_{n}+\epsilon f$, where $1_{n}$ denotes an $n\times 1$ column vector of all ones, $E$ and $f$ are given error matrix and vector generated by using the Matlab. If $B^{T}B-\sigma_{n+1}^{2}I_{n}$ is positive definite, the solution of the total least squares problem \cite{Bjorck1996} related to $B$, $d$ is as $
	x_{TLS}=(B^{T}B-\sigma_{n+1}^{2}I_{n})^{-1}B^{T}d$,
where $\sigma_{n+1}$ is the smallest singular value of $(B, d)$. It is equivalent to solving the ILS problem with $
A=\left(\begin{matrix}
	B ;
	\sigma_{n+1}I
\end{matrix}\nonumber
\right),~
b=\left(\begin{matrix}
	d;
	\mathbf{0}
\end{matrix}\nonumber
\right).$
\end{example}

\begin{figure}[htbp]
	\centering
	\includegraphics[height=6cm,width=10cm]{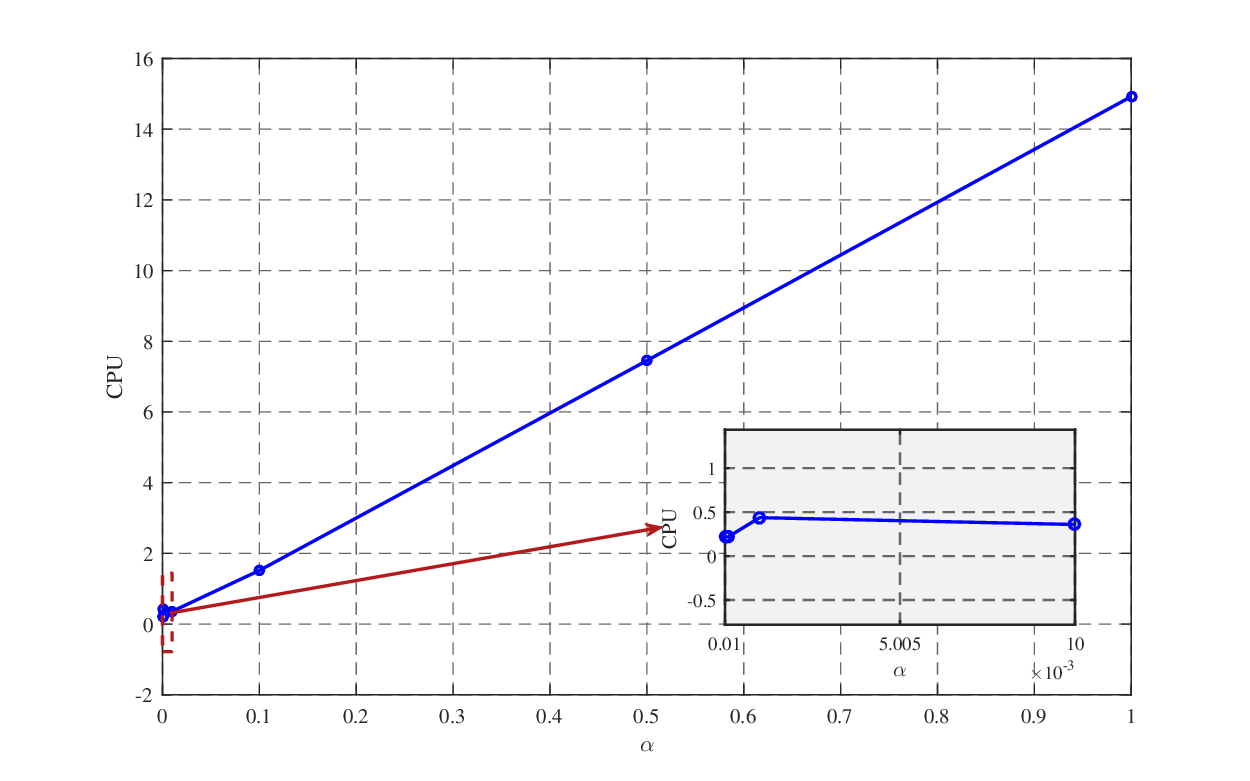}
	\caption{\footnotesize \emph{The testing parameters vs CPU times of DS iteration method for Example 4.2.}}\label{fig2}
\end{figure}

\begin{table}[htbp]  
	\centering
	\captionsetup{width=0.9\textwidth, font=small, labelfont=bf, labelsep=period}
	\renewcommand{\arraystretch}{1}  
	\setlength{\tabcolsep}{5pt}       
	\caption{The numerical results of the SP, GSP, ADI and DS iteration methods with sparse matrices for Example 4.2.}\label{tab2}
	
	\begin{tabular}{l *{8}{r}}
		\toprule
		\multirow{2}{*}{Method} 
		& \multicolumn{2}{c}{$m=1024,\; n=512$} 
		& \multicolumn{2}{c}{$m=2048,\; n=1024$} 
		& \multicolumn{2}{c}{$m=4096,\; n=2048$} 
		& \multicolumn{2}{c}{$m=8192,\; n=4096$} \\
		\cmidrule(lr){2-3} \cmidrule(lr){4-5} \cmidrule(lr){6-7} \cmidrule(lr){8-9}
		& IT & CPU (s) & IT & CPU (s) & IT & CPU (s) & IT & CPU (s) \\
		\midrule
		SP   & 1 & 0.1563 & 1 & 0.3438 & 1 & 1.1250 & 1 & 7.5313 \\
		GSP  & 1 & 0.1406 & 1 & 0.2656 & 1 & 0.9531 & 1 & 7.2188 \\
		ADI  & 1 & 0.2031 & 1 & 0.3281 & 1 & 2.1250 & 1 & 13.6563 \\
		DS   & 2 &  0.0156 & 2 &  0.1719& 2 &  0.6719 & 2 &  4.0469\\
		\bottomrule
	\end{tabular}
\end{table}

Figure \ref{fig2} illustrates the CPU runtime of the DS iteration method against the parameter  $\alpha$  in numerical experiments.
Overall, CPU time is positively correlated with $\alpha$. However, runtime grows only slowly in the small-$\alpha$ regime $\left( \alpha \in [10^{-5},10^{-2}]\right)$, as the regulatory effect of extremely small 
$\alpha$ is gently. Therefore,
 we select  $\alpha=10^{-4}$ as the optimal parameter: it lies in the stable low-runtime plateau while ensuring numerical stability, validating the high efficiency of the DS method with appropriate parameter tuning. As presented in Table \ref{tab2},  DS iteration method outperforms all compared methods in  computational efficiency for large-scale simulations.

\section{Conclusions}
In this paper, we have focused on solving the large indefinite least squares (ILS) problem. We proposed a double splitting (DS) iteration method based on splitting of the normal equation coefficient matrix $A^TJA $, which leads to a two-step iteration scheme that incorporates the previous two iterates. A detailed spectral analysis was performed to establish the convergence condition of the DS iteration amethod.  This theoretical result was further supported by a remark connecting the double splitting framework to existing convergence criteria for Hermitian positive definite matrices.

Compared with classical single splitting methods (namely SP, GSP and ADI methods), the DS iteration required only two iterations to converge under the chosen optimal parameter, while achieving significantly lower CPU times. The results confirm that the double splitting strategy can substantially accelerate the solution process for large-scale dense/sparse ILS problem.

Future work may include using comparison theorem \cite{Shen2007} to construct more effective double splitting iteration method, developing adaptive parameter selection strategies, extending the double splitting idea to other types of indefinite linear system.

\section*{Declarations}


\noindent{\bf Availability of supporting data} Data availability is not applicable to this article as no new data were created or analyzed in this study.\\

 

\noindent{\bf Funding}
This work was supported by the National Natural Science Foundation of China (No.12361082) and Doctoral Research Starting Funding of Lanzhou University of Technology (05/062506).


\end{document}